\documentclass[10pt,a4paper]{amsart}
\usepackage{latexsym}

\def\LaTeX{L\kern-.36em\raise.3ex\hbox{a}\kern-.15em
    T\kern-.1667em\lower.7ex\hbox{E}\kern-.125emX}

\def\cc{{\mathcal C}}

\def\hh{{\mathcal H}}

\def\lll{{\mathcal L}}
\def\mm{{\mathcal M}}

\def\ffi{\varphi}
\def\eps{\varepsilon}

\def\B{{\mathbb{B}}}

\def\N{{\mathbb{N}}}
\def\P{{\mathbb{P}}}

\def\R{{\mathbb{R}}}
\def\S{{\mathbb{S}}}

\def\X{{\mathbb{X}}}
\def\Y{{\mathbb{Y}}}

\newcommand{\norm}[1]{{\left\|{#1}\right\|}}
\newcommand{\ent}[1]{{\left[{#1}\right]}}
\newcommand{\abs}[1]{{\left|{#1}\right|}}

\newenvironment{theorem}[1][]{\ \\ \vskip3pt\noindent\sl\textbf{Theorem\ #1}\ }{\rm}

\newenvironment{definition}[1][]{\ \\ \vskip3pt\noindent\sl\textbf{Definition}\
}{\rm\vskip3pt}

\renewenvironment{proof}[1][]{\noindent\textit{#1}.\ \rm}{}
\newenvironment{proposition}[1][]{\ \\ \vskip3pt\noindent\sl\textbf{Proposition\ #1}\
}{\rm\vskip3pt}

\newenvironment{lemmenum}[1][]{\ \\ \vskip1pt\noindent\sl\textbf{Lemma\ #1}\
}{\rm\vskip3pt}

\newenvironment{notation}[1][]{\ \\ \vskip1pt\noindent\rm\textit{Notation}\ :\
}{\rm\vskip2pt}

\newenvironment{remarque}[1][]{\ \\ \vskip1pt\noindent\rm\textit{Remark}\ :\
}{\rm\vskip2pt}
\newenvironment{remarquenum}[1][]{\ \\ \vskip1pt\noindent\rm\textit{Remark\ #1}\ :\
}{\rm\vskip2pt}

\newenvironment{corollaire}[1][]{\ \\ \vskip3pt\noindent\sl\textbf{Corollary\ #1}\
}{\rm\vskip3pt}

\begin{document}
\title[Harmonic functions on the real hyperbolic ball]{Harmonic functions on the real hyperbolic ball I
:\\ 
Boundary values and atomic decomposition of Hardy spaces}
\author{Philippe JAMING$^1$}
\address{Universit\'e d'Orl\'eans\\ Facult\'e des Sciences\\ 
D\'epartement de Math\'ematiques\\BP 6759\\ F 45067 ORLEANS Cedex 2\\
FRANCE}
\email{jaming@labomath.univ-orleans.fr}
\date{}
\thanks{$^1$Most of the results in this paper are part of my Ph.D.
thesis ``Trois probl\`emes d'analyse harmonique'' written at the
University of Orl\'eans under the direction of Aline Bonami, to whom I
wish to express my sincere gratitude. I also want to thank
Sandrine Grellier for valuable conversations.}

\begin{abstract} In this article we study harmonic functions for the
Laplace-Beltrami
operator on the real hyperbolic space $\B_n$. We obtain necessary and
sufficient conditions for this functions and their normal derivatives 
to have a boundary distribution. In doing so, we put forward different
behaviors of
hyperbolic harmonic functions according to the parity of the dimension
of the hyperbolic ball $\B_n$. We then study Hardy spaces $H^p(\B_n)$,
$0<p<\infty$, whose elements appear as the hyperbolic harmonic
extensions of distributions belonging to the
Hardy spaces of the sphere $H^p(\S^{n-1})$. In particular, we obtain an
atomic decomposition of this spaces.
\end{abstract}

\maketitle

\section{Introduction}

In this article, we study boundary behavior of harmonic functions on
the real hyperbolic ball, partly in view of establishing a theory of
Hardy and Hardy-Sobolev spaces of such functions.

While studying Hardy spaces of Euclidean harmonic functions on the unit
ball $\B_n$ of $\R^n$, one is often lead to consider estimates of this
functions on balls with radius smaller than the distance of the center
of the ball to the boundary $\S^{n-1}$ of $\B_n$. Thus hyperbolic
geometry is implicitly used for the study of Euclidean harmonic
functions, in particular when one considers boundary behavior. As Hardy
spaces of Euclidean harmonic functions are the spaces of Euclidean
harmonic extensions of distributions in the Hardy spaces on the sphere,
it is tempting to study these last spaces directly through their
hyperbolic harmonic extension.

The other origin of this paper is the study of Hardy and Hardy-Sobolev
spaces of $\mm$-harmonic functions related to the complex hyperbolic
metric on the unit ball, as exposed in \cite{ABC} and \cite{BBG}. Our
aim is to develop a similar theory in the case of the real hyperbolic
ball. In the sequel, $n$ will be an integer, $n\geq3$ and $p$ a real
number, $0<p<\infty$.

Let $SO(n,1)$ be the Lorenz group. It is well known that $SO(n,1)$ acts
conformly on $\B_n$. The corresponding
Laplace-Beltrami operator, invariant for the considered action, is given
by

$$D=(1-\abs{x}^2)^2\Delta+2(n-2)(1-\abs{x}^2)N$$
with $\Delta$ the Euclidean laplacian and
$N=\sum_{i=1}^nx_i\frac{\partial}{\partial x_i}$ the normal derivation
operator. Functions $u$ that are harmonic for this laplacian will be
called $\hh$-harmonic. The ``hyperbolic'' Poisson kernel that solves the
Dirichlet problem for $D$ is defined for $x\in\B_n$ and
$\xi\in\S^{n-1}$ by

$$\P_h(x,\xi)=\left(\frac{1-\abs{x}^2}{1+\abs{x}^2-2<x,\xi>}\right)^{n-1}.$$

With help of this kernel, one can extend distributions on $\S^{n-1}$ to
$\hh$-harmonic functions on $\B_n$ in the same way as the Euclidean
Poisson kernel extends distributions on $\S^{n-1}$ to Euclidean harmonic
functions on $\B_n$. Our first concern is to determine which
$\hh$-harmonic functions are obtained in this way. We then study the
boundary behavior of their normal derivatives. In doing so, we put
forward that, in odd dimension, normal derivatives of $\hh$-harmonic
functions behave similarly to $\mm$-harmonic functions whereas they
behave like Euclidean harmonic functions in even dimension.

Finally, define $H^p(\S^{n-1})$ as $L^p(\S^{n-1})$ if $1<p<\infty$ and
as the real analog of Garnett-Latter's atomic $H^p$ space if $p\leq1$. Let
$H^p(\B_n)$ be the space of Euclidean harmonic functions $\B_n$ such
that $\zeta\mapsto\sup_{0<r<1}\abs{u(r\zeta)}\in L^p(\S^{n-1})$.
Garnett-Latter's theorem asserts that this space is the space of
Euclidean harmonic extensions of distributions in $H^p(\S^{n-1})$. We
prove here that the space $\hh^p(\B_n)$ of $\hh$-harmonic functions
such that $\zeta\mapsto\sup_{0<r<1}\abs{u(r\zeta)}\in L^p(\S^{n-1})$ is
the space of $\hh$-harmonic extensions of distributions in
$H^p(\S^{n-1})$.
\vskip6pt
This article is organized as follows : in section 2 we present the
setting of the problem and a few preliminary results. Section 3 is
devoted to the study of boundary behavior of $\hh$-harmonic functions
and concludes with the study of the behavior of their normal
derivatives. We conclude in section 4 with the atomic decomposition
theorem.

\section{Setting}

\subsection{$SO(n,1)$ and its action on $\B_n$}
Let $SO(n,1)\subset GL_{n+1}(\R)$, ($n\geq3$) be the identity component
of the group of matrices
$g=(g_{ij})_{0\leq i,j\leq n}$ such that $g_{00}\geq1$, $\det g=1$ and
that leaves invariant the quadratic form $-x_0^2+x_1^2+\ldots+x_n^2$.

Let $\abs{.}$ be the Euclidean norm on $\R^n$, $\B_n=\{x\in\R^n\
:\ \abs{x}<1\}$ and $\S^{n-1}=\partial\B_n=\{x\in\R^n\ :\ \abs{x}=1\}$.
It is well known ({\it cf.} \cite{Sam}) that $SO(n,1)$ acts conformaly
on $\B_n$. The action is given by $y=g.x$ with

$$y_p=\frac{\frac{1+\abs{x}^2}{2}g_{p0}+\sum_{l=1}^ng_{pl}x_l}{\frac{1-\abs{x}^2}{2}+
\frac{1+\abs{x}^2}{2}g_{00}+\sum_{l=1}^ng_{0l}x_l}
\qquad\mathrm{for\ }p=1,\ldots,n.$$
The invariant measure on $\B_n$ is given by

$$d\mu=\frac{dx}{(1-\abs{x}^2)^{n-1}}=
\frac{r^{n-1}drd\sigma}{(1-r^2)^{n-1}}$$
where $dx$ is the Lebesgue measure on $\B_n$ and $d\sigma$ is the
surface measure on $\S^{n-1}$.

We will need the following fact about this action (see \cite{Jathese}):

{\bf Fact 1} {\sl Let $g\in SO(n,1)$ and let $x_0=g.0$. If
$0<\eps<\frac{1}{6}$, then}

$$B\bigl(x_0,\frac{\sqrt{2}}{8}(1-\abs{x_0}^2)\eps\bigr)\subset
g.B(0,\eps)\subset B\bigl(x_0,6(1-\abs{x_0}^2)\eps\bigr).$$

\subsection{The invariant laplacian on $\B_n$ and the associated Poisson
kernel}
%\markboth{Fonctions $\hh$-harmoniques}{Laplacien invariant sur $\B_n$}

From \cite{Sam} we know that the invariant laplacian on $\B_n$
for the considered action can be written as
%\index{Laplacien Invariant sur $\B_n$, $D$}

$$D=(1-r^2)^2\Delta+2(n-2)(1-r^2)\sum_{i=1}^nx_i\frac{\partial}{\partial
x_i}$$%\index{daaa@$D$}
where $r=\abs{x}=(x_1^2+\ldots+x_n^2)^{1/2}$ and $\Delta$
is the Euclidean laplacian
$\Delta=\sum_{i=1}^n\frac{\partial^2}{\partial x_i^2}$.

Note that $D$ is given in radial-tangential coordinates by

$$D=\frac{1-r^2}{r^2}\ent{(1-r^2)N^2+(n-2)(1+r^2)N+(1-r^2)\Delta_\sigma}$$
with $N=r\frac{d}{dr}=\sum_{i=1}^nx_i\frac{\partial}{\partial x_i}$ and $\Delta_\sigma$ the tangential part of the
Euclidean laplacian.

\begin{definition} A function $u$ on $\B_n$ is
$\hh$-harmonic if $Du=0$ on $\B_n$.
\end{definition}

The Poisson kernel that solves the Dirichlet problem  associated to $D$ is

$$\P_h(r\eta,\xi)=\left(\frac{1-r^2}{1+r^2-2r<\eta,\xi>}\right)^{n-1}$$
for $0\leq r<1$, $\eta,\xi\in\S^{n-1}$ {\it i.e.} for
$r\eta\in\B_n$ and $\xi\in\S^{n-1}$.
%\index{phreta@$\P_h(r\eta,\xi)$}

Recall that the Euclidean Poisson kernel on the ball is given by

$$\P_e(r\eta,\xi)=\frac{1-r^2}{(1+r^2-2r<\eta,\xi>)^{\frac{n}{2}}}$$
%\index{pereta@$\P_e(r\eta,\xi)$}

\begin{notation} For a distribution $\ffi$ on $\S^{n-1}$, we define
$\P_e\ent{\ffi}:\B_n\mapsto\R$ and $\P_h\ent{\ffi}:\B_n\mapsto\R$ by

\begin{align}
\P_e\ent{\ffi}(r\eta)=&<\ffi,\P_e(r\eta,.)>\notag\\
\P_h\ent{\ffi}(r\eta)=&<\ffi,\P_h(r\eta,.)>\notag
\end{align}
$\P_e\ent{\ffi}$ is the {\it Poisson integral} of $\ffi$, and
$\P_h\ent{\ffi}$ will be called the {\it $\hh$-Poisson integral} of
$\ffi$.
\end{notation}

Finally, $\hh$-harmonic functions satisfy mean value equalities :
{\sl let $a\in\B_n$ and $g\in SO(n,1)$ such that $g.0=a$. Then, for 
every $\hh$-harmonic function $u$,}

$$u(a)=\frac{1}{\mu\bigl(B(0,r)\bigr)}
\int_{g.B(0,r)}u(x)d\mu(x).$$
Thus, with fact 1 and $d\mu=\frac{dx}{(1-\abs{x}^2)^{n-1}}$, we see
that

$$\abs{u(a)}\leq\frac{C}{(1-\abs{a}^2)^n}
\int_{B\bigl(a,6(1-\abs{a}^2)\eps\bigr)}\abs{u(x)}dx
\eqno(2.1)$$

\subsection{Expansion of $\hh$-harmonic functions in spherical
harmonics}
%\markboth{Fonctions $\hh$-harmoniques}{D\'eveloppement en harmoniques sph\'eriques}

\begin{notation} For $a\in\R$, write
$(a)_k=\frac{\Gamma(a+k)}{\Gamma(a)}$ thus $(a)_0=1$ and 
$(a)_k=a(a+1)\ldots(a+k-1)$ if $k=1,2,\ldots$. For
$a,b,c$ three real parameters, ${}_2F_1$ denotes Gauss' {\it
hyper-geometric function} defined by

$${}_2F_1(a,b,c;x)=\sum_{k=0}^\infty\frac{(a)_k(b)_k}{(c)_k k!}x^k.$$

Let $F_l(x)={}_2F_1(l,1-\frac{n}{2},l+\frac{n}{2};x)$ and
$f_l(x)=\frac{F_l(x)}{F_l(1)}$. (See \cite{Erd} for properties of
${}_2F_1$).
\end{notation}

\begin{remarque} If $n>2$ is even, $1-\frac{n}{2}$ is a negative integer
thus ${}_2F_1(l,1-\frac{n}{2},l+\frac{n}{2},r^2)$ is a polynomial
in $r$ of degree $n$.
\end{remarque}

In \cite{Mi1}, \cite{Mi2} and \cite{Sam}, the spherical harmonic expansion
of $\hh$-harmonic functions has been obtained. An other proof based on
\cite{ABC} can be found
in \cite{Jathese}. We have the following :

\begin{theorem}[1] Let $u$ be an $\hh$-harmonic function of class
$\cc^2$ on $\B_n$. Then the spherical harmonic expansion of $u$ is given
by

$$u(r\zeta)=\sum_lF_l(r^2)u_l(r\zeta),$$
where this series is absolutely convergent and uniformly convergent
on every compact subset of $\B_n$.

Moreover if $\ffi\in\cc(\S^{n-1})$, the Dirichlet
problem
$\left\{\begin{array}{cc}
Du=0&\mathrm{\ in\ }\B_n\\
u=\ffi&\mathrm{\ on\ }\S^{n-1}\\
\end{array}\right.$
has a unique solution $u\in\cc(\overline{\B_n})$ given by

$$u(z)=\int_{\S^{n-1}}\ffi(\zeta)\P_h(z,\zeta)d\sigma(\zeta)=\P_h\ent{\ffi}(z)$$
also given by

$$u(r\zeta)=\sum_l f_l(r^2)r^l\ffi_l(\zeta)$$
where $\ffi=\sum_l\ffi_l$ is the spherical harmonic expansion of $\ffi$.
\end{theorem}

\section{Boundary values of $\hh$-harmonic functions}
\label{distribord}

{\sl In this chapter we prove results about the behavior on the 
boundary of $\hh$-harmonic functions and their normal derivatives. For
$\hh$-harmonic functions, the results are similar to the results for
Euclidean harmonic
functions. On the opposite, for the normal derivatives of  $\hh$-harmonic
functions, the boundary behavior depends on the dimension of the
space.}

\subsection{Definition of Hardy spaces}

\begin{notation} For $u$ a function defined on $\B_n$, define the
 {\it radial maximal function} 
$\mm\ent{u}:\S^{n-1}\mapsto\R_+$ by

$$\mm\ent{u}(\zeta)=\sup_{0<t<1}\abs{u(t\zeta)}.$$
\end{notation}

We will now study $\hh^p$ spaces of $\hh$-harmonic functions
defined as follows :

\begin{definition} Let $0<p<\infty$. Let $\hh^p$ be the space of
$\hh$-harmonic functions $u$ such that $\mm\ent{u}\in L^p(\S^{n-1})$,
endowed with the ``norm'' 

$$\norm{u}_{\hh^p}=\norm{\mm
u}_{L^p(\S^{n-1})}=\norm{\sup_{0<t<1}\abs{u(t.)}}_{L^p(\S^{n-1})}.$$
We will call $\hh^p$ the Hardy space of $\hh$-harmonic functions.
\end{definition}
%\index{hp@$\hh^p$}
%\index{normhp@$\norm{.}_{\hh^p}$}

\begin{remarque} If $0<p<1$, the application $u\mapsto\norm{u}_{\hh^p}$
is not a norm, however the application
$u,v\mapsto\norm{u-v}_{\hh^p}$ defines a metric on $\hh^p$.
In the sequel, we will often use the abuse of language to call
$\norm{.}_{\hh^p}$ a norm whatever $p$ might be.
\end{remarque}

\begin{definition} A function $u$ on $\B_n$ is said to have a
distribution boundary value if for every $\Phi\in\cc^\infty(\S^{n-1})$,
the limit

$$\lim_{r\rightarrow1}\int_{\S^{n-1}}u(r\zeta)\Phi(\zeta)d\sigma(\zeta)$$
exists. In case $u$ is $\hh$-harmonic, this is equivalent to the
existence of a distribution $f$ such that $u=\P_h\ent{f}$.
\end{definition}

\subsection{Boundary distributions of functions in $\hh^p$}
\label{rangun}

In this section, we are going to characterize boundary values of
functions in
$\hh^p$. The characterizations we obtain are similar to those obtained
for
harmonic functions on $\R^{n+1}_+$ or for $\mm$-harmonic functions. The
proofs
are inspired by \cite{ABC} and \cite{FS1}.

The first result concerns functions in $\hh^p$, $p\geq 1$.

\begin{proposition}[2] Let $u$ be an $\hh$-harmonic function.

\begin{enumerate}
\item If $1<p<\infty$, then
$$\sup_{0<r<1}\int_{\S^{n-1}}\abs{u(r\zeta)}^pd\sigma(\zeta)<+\infty$$
if and only if there exists $f\in L^p(\S^{n-1})$ such that 
$u=\P_h\ent{f}$.

\item For $p=1$,
$$\sup_{0<r<1}\int_{\S^{n-1}}\abs{u(r\zeta)}d\sigma(\zeta)<+\infty$$
if and only if there exists a measure $\mu$ on $\S^{n-1}$ such that
$u=\P_h\ent{\mu}$.
\end{enumerate}
\end{proposition}

\begin{proof}[Proof] Assume that $u=\P_h\ent{f}$ with $f\in
L^p(\S^{n-1})$. As

$$\norm{\P_h(r\zeta,.)}_{L^1(\S^{n-1})}=1,$$
H\"older's inequality gives

$$\abs{u(r\zeta)}^p\leq\int_{\S^{n-1}}\P_h(r\zeta,\xi)\abs{f(\xi)}^p
d\sigma(\xi)=\int_{\S^{n-1}}\P_h(\zeta,r\xi)\abs{f(\xi)}^p
d\sigma(\xi)$$
an integration in $\zeta$ and Fubini leads to the desired result.

Conversely, if the $L^p(\S^{n-1})$ norms of $\zeta\mapsto u(r\zeta)$
are uniformly bounded, there exists a sequence $r_m\rightarrow 1$ and a
function $\ffi\in L^p$ such that $u(r_m\zeta)\rightarrow \ffi(\zeta)$
$*$-weakly thus weakly in $L^p(\S^{n-1})$. But then, for
$r\zeta\in\B_n$ fixed,

\begin{align}
\P_h\ent{\ffi}(r\zeta)=&\lim_{m\rightarrow+\infty}
\int_{\S^{n-1}}\P_h(r\zeta,\xi)u(r_m\xi)d\sigma(\xi)\notag\\
=&\lim_{m\rightarrow+\infty}\sum_{l\geq0}\frac{F_l(r_m^2)}{F_l(1)}r_m^l
\int_{\S^{n-1}}\P_h(r\zeta,\xi)u_l(\xi)d\sigma(\xi)\notag\\
=&\lim_{m\rightarrow+\infty}\sum_{l\geq0}\frac{F_l(r_m^2)}{F_l(1)}
r_m^lf_l(r)r^lu_l(\zeta)\notag\\
=&\sum_{l\geq0}f_l(r)r^lu_l(\zeta)=u(r\zeta).\notag
\end{align}

The proof in the case $p=1$ is obtained in a similar fashion
using the duality $\bigl(L^1,\mm(\S^{n-1})\bigr)$.\hfill$\Box$
\end{proof}

We are now going to prove that $\hh$-harmonic functions have a boundary
distribution if and only if they satisfy a given growth condition. For
this, we will need the folowing lemma ( \cite{ABC}, lemma 10).

\begin{lemmenum}[3] Let $F\in\cc^2\left(\ent{\frac{1}{2},1}\right)$
and
$h\in\cc^1\left(\ent{\frac{1}{2},1}\right)$. Assume that

$$F''(x)+\frac{h(x)}{1-x}F'(x)=O(1-x)^{-\alpha}$$
when $x\rightarrow1$. Then
\begin{enumerate}
\item If $\alpha >2$ then $F(x)=O(1-x)^{-\alpha+1}$.

\item If $1<\alpha<2$ Then $\lim_{x\rightarrow1}F(x)$ exists.
\end{enumerate}
\end{lemmenum}

We are now in position to prove

\begin{theorem}[4] Let $u$ be an $\hh$-harmonic function. Then
$u$
admits a boundary value in the sense of distributions if and only if
there exists a constant $A$ such that

$$u(r\zeta)=O(1-r)^{-A}.$$ 
\end{theorem}

\begin{proof}[Proof] Recall that

$$D=\frac{1-r^2}{r^2}\ent{(1-r^2)N^2+(n-2)(1+r^2)N+
(1-r^2)\Delta_\sigma}\eqno(3.1)$$

Assume that $Du=0$ and that $u(r\zeta)=O\bigl((1-r)^{-A}\bigr)$. Let 
$\ffi\in\cc^{\infty}(\S^{n-1})$ and let

$$F(r)=\int_{\S^{n-1}}u(r\zeta)\ffi(\zeta)d\sigma(\zeta).$$
Formula (3.1) with $Du=0$ tells us that

$$(1-r^2)N^2F+(n-2)(1+r^2)NF+(1-r^2)\Delta_\sigma F=0.$$
where $\Delta_\sigma F$ stands for
$$\Delta_\sigma F(r)=\int_{\S^{n-1}}\Delta_\sigma
u(r\zeta)\ffi(\zeta)d\sigma(\zeta)
=\int_{\S^{n-1}}u(r\zeta)\Delta_\sigma^*\ffi(\zeta)d\sigma(\zeta)$$
with $\Delta_\sigma^*$ the adjoint operator to $\Delta_\sigma$.
Recall that $N=r\frac{d}{dr}$ thus

$$r^2F''(r)+\frac{(n-1)+(n-3)r^2}{1-r^2}rF'(r)+\Delta_\sigma F=0
\eqno(3.2)$$

Write $\psi=-\Delta_\sigma^*\ffi$ and $T$ the differential operator

$$T=r^2\frac{d^2}{dr^2}+\frac{(n-1)+(n-3)r^2}{1-r^2}r\frac{d}{dr}$$
so that equation (3.2) reads

$$TF(r)=\int_{\S^{n-1}}u(r\zeta)\psi(\zeta)d\sigma(\zeta).$$
One then immediately deduces the existence for $k=1,2,\ldots$
of a function $\psi_k\in\cc^\infty(\S^{n-1})$ such that

$$T^kF(r)=\int_{\S^{n-1}}u(r\zeta)\psi_k(\zeta)d\sigma(\zeta).$$
But we assumed that
$u(r\zeta)=O(1-r)^{-A}$. We thus have

$$T^kF(r)=O(1-r)^{-A}$$
and applying lemma 3 we obtain

$$T^{k-1}F(r)=O(1-r)^{-A+1}.$$
Therefore, starting from $T^k$ with $k=\ent{A}+1$ and
iterating the
process $k$ times, one gets that $\lim_{r\rightarrow 1}F(r)$
exists.

Conversely, if $u$ admits a boundary distribution $f$, then 
$u=\P_h\ent{f}$ {\it i.e.} $u(r\zeta)=<f,\P_h(r\zeta,.)>$. But then
$f$ being a compactly supported distribution, it is of finite order,
thus there exists $k\geq0$ such that

$$\abs{u(r\zeta)}=\abs{<f,\P_h(r\zeta,.)>}\leq
C\norm{\nabla^k_\xi\P_h(r\zeta,.)}_{L^\infty}
\leq\frac{C}{(1-r)^{n-1+k}}$$
which gives the desired estimate.\hfill$\Box$
\end{proof}

\begin{proposition}[5] Let $0<p<+\infty$ and $u$ be an 
$\hh$-harmonic function. Assume that

$$\sup_{0<r<1}\int_{\S^{n-1}}\abs{u(r\zeta)}^pd\sigma(\zeta)<\infty.$$
Then, there exists a constant $C$ such that for every $a\in\B_n$,

$$\abs{u(a)}\leq\frac{C}{(1-\abs{a})^{\frac{n-1}{p}}}.$$
In particular, $u$ has a boundary distribution $f$ {\it i.e.}
$u=\P_h\ent{f}$.
\end{proposition}

\begin{proof}[Proof] The mean value inequality implies that

$$\abs{u(a)}^p\leq\frac{C}{(1-\abs{a})^n}
\int_{B\bigl(a,(1-\abs{a})\eps\bigr)}\abs{u(x)}^pdx$$
for $\eps$ small enough. But
$B\bigl(a,(1-\abs{a})\eps\bigr)\subset\{r\zeta\ :\
(1-\eps)(1-\abs{a})\leq1-r\leq(1+\eps)(1-\abs{a})\}$ thus

$$\abs{u(a)}^p\leq\frac{C}{(1-\abs{a})^n}
\int_{1-(1+\eps)(1-\abs{a})}^{1-(1-\eps)(1-\abs{a})}\int_{\S^{n-1}}
\abs{u(r\zeta)}^pd\sigma(\zeta)r^{n-1}dr
\leq\frac{C}{(1-\abs{a})^{n-1}}.\eqno\Box$$
\end{proof}

\begin{remarque} Theorem 4 is well known. It has been proved 
by J.B. Lewis \cite{Lew} in the case of symmetric spaces of rank 1 and
eigenvectors of the Laplace-Beltrami operator (for arbitrary eigenvalues) and further
generalized by E.P. van den Ban and H. Schlichtkrull \cite{vdBS}.
\end{remarque}

\subsection{Distribution boundary values of $\hh$-harmonic functions.}
%\markboth{Fonctions $\hh$-harmoniques}{Distributions au bord}

\begin{notation} For $1\leq i,j\leq n$, $i\not=j$, let
$\lll_{i,j}=x_i\frac{\partial}{\partial x_j}-
x_j\frac{\partial}{\partial x_i}$. Then the $\lll_{i,j}$'s commute and
commute with $N$. Further, if $u$ is $\hh$-harmonic, then $\lll_{i,j}$
is also $\hh$-harmonic. Finally, $N$ and $\{\lll_{i,j}\}_{1\leq
i\not=j\leq n}$ generate $\nabla^k$ outside a neighbourhood of the
origin.
\end{notation}

Recall that $Du=0$ if and only if

$$(1-r^2)N^2u+(n-2)(1+r^2)Nu+(1-r^2)\Delta_\sigma u=0.
\eqno(3.3)$$
Apply $N^{k-1}$ on both sides of this equality and isolate terms of order
$k+1$ and $k$ :

\begin{align}
(1-r^2)N^{k+1}u-2(k-1)r^2N^ku+&(n-2)(1+r^2)N^ku\notag\\
=&r^2\sum_{j=0}^{k-3}\begin{pmatrix}k-1\\j\\ \end{pmatrix}2^{k-j-1}N^{j+2}u
+r^2\sum_{j=0}^{k-2}\begin{pmatrix}k-1\\j\\ \end{pmatrix}2^{k-j-1}N^{j}
\Delta_\sigma u\notag\\
&-(n-2)r^2\sum_{j=0}^{k-2}\begin{pmatrix}k-1\\j\\ \end{pmatrix}2^{k-j-1}N^{j+1}u
-(1-r^2)N^{k-1}\Delta_\sigma u\tag{3.4}
\end{align}

We are now in position to prove the following lemma :

\begin{lemmenum}[6] Let $u$ be an $\hh$-harmonic function with a
boundary distribution. Let $\Y$ be a product of operators of the form
$\lll_{i,j}$ and let $\X=N^k\Y$. Then if $k\leq n-2$, $\X u$ has a 
distribution boundary value in the sense that

$$\lim_{r\rightarrow 1}\int_{\S^{n-1}}\X
u(r\zeta)\Phi(\zeta)d\sigma(\zeta)$$
exists for every function $\Phi\in\cc^\infty(\S^{n-1})$.

If $k=n-1$, the previous integral is a
$O\left(\log\frac{1}{1-r}\right)$, in 
particular

$$\lim_{r\rightarrow 1}(1-r^2)
\int_{\S^{n-1}}\X u(r\zeta)\Phi(\zeta)d\sigma(\zeta)=0.$$
\end{lemmenum}

\begin{remarquenum}[1] If $u$ has a 
boundary distribution, then $\lll_{i,j} u$ has a boundary distribution.
\end{remarquenum}

\begin{remarquenum}[2] As $\nabla^k$ is generated outside a
neighbourhood of the origin by operators of the form
$N^l\Y$
where $\Y$ is a product of at most $k-l$ operators of the form
$\lll_{i,j}$,
We deduce from the lemma that if
$k\leq n-2$, $\nabla^k$ has a boundary distribution, whereas

$$\int_{\S^{n-1}}\nabla^{n-1}u(r\zeta)\Phi(\zeta)d\sigma(\zeta)$$
has {\sl a priori} logarithmic growth.
\end{remarquenum}

\begin{proof}[Proof] Proceed by induction on $k$. Fix
$\Phi\in\cc^\infty(\S^{n-1})$ and let $\Y$ be a product of operators of
the
form $\lll_{i,j}$. Let

$$\psi_k(r)=\int_{\S^{n-1}}N^k\Y
u(r\zeta)\Phi(\zeta)d\sigma(\zeta);\qquad0<r<1.$$

Applying $\Y$ to formula (3.4) and noticing that 
$\Y$ and $N$ commute, the induction hypothesis implies that the function

$$g(r)=(1-r^2)N\psi_k(r)-2(k-1)r^2\psi_k(r)+(n-2)(1+r^2)\psi_k(r)
\eqno(3.5)$$
has a limit $L$ when $r\rightarrow1$.

But, solving the differential equation (3.5),
($N=r\frac{d}{dr}$), we get

$$\psi_k(r)=\lambda\frac{(1-r^2)^{n-k-1}}{r^{n-2}}+
\frac{1}{r^{n-2}}(1-r^2)^{n-k-1}
\int_0^r\frac{g(s)s^{n-3}}{(1+s)^{n-k}}(1-s)^{-(n-k-1)-1}ds.$$
Thus, if $k<n-1$, we obtain that $\psi_k(r)$ has limit
$\frac{L}{2(n-k-1)}$
whereas if $k=n-1$, $\psi_k(r)$ has logarithmic growth.\hfill$\Box$
\end{proof}

\begin{remarque} We will show at the end of this section that if $n$ is
even,
$N^{n-1}u$ can have a better than logarithmic growth, whereas if $n$ is
odd, only constant functions have a better than logarithmic growth.
\end{remarque}

\begin{corollaire}[7] Let $P_k$ be the sequence of polynomials
defined
by $P_0=2(n-1)$, $P_1=0$ and for $2\leq k\leq n$,

\begin{align}
P_k(X)=&2^{k-1}(k-1)!\sum_{j=2}^{k-2}
\frac{n(j-1)-(n-2)k}{2^j(n-j-1)(k-j+1)!(j-1)!}P_j(X)\notag\\
&+2^{k-2}(k-1)!\sum_{j=2}^{k-3}\frac{1}{2^j(n-j-1)(k-j-1)!j!}XP_j(X)
+2^{k-1}X\notag
\end{align}
Then, for every $\hh$-harmonic function $u$ having a distribution
boundary value, and for every $1\leq k\leq n-2$,
$N^ku=\frac{1}{2(n-k-1)}P_k(\Delta_\sigma)u$ as boundary distributions, 
{\it i.e.} for every $\Phi\in\cc^\infty(\S^{n-1})$,

$$\lim_{r\rightarrow1}
\int_{\S^{n-1}}\left(N^ku(r\zeta)-\frac{1}{2(n-k-1)}P_k(\Delta_\sigma)
u(r\zeta)\right)\Phi(\zeta)d\sigma(\zeta)=0.$$
\end{corollaire}

\begin{proof}[Proof] For convenience, write $Q_k=\frac{1}{2(n-k-1)}P_k$. 
As $n\geq 3$, for $u$ $\hh$-harmonic having a boundary distribution, formula 
(3.3) and lemma 6 impliy that $Nu=0$
on the boundary, thus the result for $k=1$.

Next, notice that $N^ku=Q_k(\Delta_\sigma)u$ on the boundary implies
$\Delta_\sigma N^ku=\Delta_\sigma Q_k(\Delta_\sigma)u$ on the boundary.

Assume now that $N^ju=Q_j(\Delta_\sigma)u$ on the boundary for $j\leq k-1$.
If $k\leq n-2$, lemma 6 tells us that $(1-r^2)N^{k+1}u=0$ on the boundary and
that $(1-r^2)N^{k-1}\Delta_\sigma u=0$ on the boundary. Formula (3.4) gives 
then, when $r\rightarrow1$,

\begin{align}
\bigl(-2(k-1)+2(n-2)\bigr)N^ku=&\sum_{j=0}^{k-3}\begin{pmatrix}k-1\\j\\
\end{pmatrix}2^{k-j-1}N^{j+2}u+\sum_{j=0}^{k-3}\begin{pmatrix}k-1\\j\\
\end{pmatrix}2^{k-j-1}N^j\Delta_\sigma
u\notag\\
&-(n-2)\sum_{j=0}^{k-2}\begin{pmatrix}k-1\\j\\ \end{pmatrix}
2^{k-j-1}N^{j+1}u.\notag
\end{align}
But, by the induction hypothesis,
$N^ju=Q_j(\Delta_\sigma)u$ and with the previous remark
$N^j\Delta_\sigma u=\Delta_\sigma N^ju=\Delta_\sigma
Q_j(\Delta_\sigma)u$, therefore

\begin{align}
\bigl(-2(k-1)+2(n-2)\bigr)N^ku=&\sum_{j=0}^{k-3}\begin{pmatrix}k-1\\j\\
\end{pmatrix}2^{k-j-1}Q_{j+2}(\Delta_\sigma)u+
\sum_{j=0}^{k-3}\begin{pmatrix}k-1\\j\\ \end{pmatrix}
2^{k-j-1}\Delta_\sigma Q_j(\Delta_\sigma)u\notag\\
&-(n-2)\sum_{j=0}^{k-2}\begin{pmatrix}k-1\\j\\ \end{pmatrix}
2^{k-j-1}Q_{j+1}(\Delta_\sigma)u.\notag
\end{align}
finally, using $Q_0=1$ and $Q_1=0$ and grouping terms, we get
the desired result.\hfill$\Box$
\end{proof}

\begin{remarquenum}[1] One easily sees that $P_k$ is a polynomial of degree
$\ent{\frac{k}{2}}$ and that for $k\geq2$, $P_k$ has no constant term.
\end{remarquenum}

\begin{remarquenum}[2] According to corollary 7, $Nu=0$ on the
boundary.
On the other hand, an easy computation leads to $DNu=-4(n-2)Nu$ {\it i.e.}
$Nu$ is an eigenvector of $D$ for an eigenvalue of the form $(s^2-1)(n-1)^2$
(with $s=\frac{n-3}{n-1}$) thus $(s+1)\frac{n-1}{2}=n-2\in\N^*$. This is
precisely the case where it is impossible to reconstruct $Nu$ with help of a
convolution by a power of the Poisson kernel (see \cite{Sam}).
\end{remarquenum}

\begin{remarquenum}[3] The fact that for every $\hh$-harmonic function
$u$, 
$Nu=0$ on the boundary is in strong contrast with Euclidean harmonic
functions.
Actually, if $v$ is an Euclidean harmonic function on $\B_n$, and if
$Nv=0$ 
on the boundary, then $v$ is a constant.
\end{remarquenum}

\subsection{Boundary distribution of the $n-1^{th}$ derivative}
%\markboth{Fonctions $\hh$-harmoniques}{Distributions au bord des d\'eriv\'ees $n-1$-i\`emes}

{\sl In this section we prove that, in odd dimension, normal derivatives
of $\hh$-harmonic functions have a boundary behavior similar to the
complex case of $\mm$-harmonic functions as exhibited in \cite{BBG}
(with pluriharmonic functions playing the role of constant functions)
whereas, in even dimension, the behavior is similar to the Euclidean
harmonic case.}

\begin{theorem}[8] $\diamond$ Assume $n$ is odd.

Let $u$ be an $\hh$-harmonic function having a boundary distribution.
The 
following assertions are equivalent :
\begin{enumerate}
\item $u$ is a constant,
\item $N^{n-1}u$ has a boundary distribution,
\item $\int_{\S^{n-1}}N^{n-1}u(r\zeta)\Phi(\zeta)d\sigma(\zeta)=
o\left(\log\frac{1}{1-r}\right)$ for every
$\Phi\in\cc^\infty(\S^{n-1})$.

$\diamond$ Assume $n$ is even, then if $\ffi\in\cc^\infty(\S^{n-1})$,
$\P_h\ent{\ffi}\in\cc^\infty(\overline{\B_n})$. In particular, if $u$ is
$\hh$-harmonic with a boundary distribution, then for every $k\geq0$, $N^ku$ has
a boundary distribution.
\end{enumerate}
\end{theorem}

\begin{proof}[Proof] $\diamond$ Assume first $n$ is odd. The implications 
$(1)\Rightarrow(2)$ and $(2)\Rightarrow(3)$ being obvious, let us prove
$(3)\Rightarrow(1)$. Theorem 1 tells us that an $\hh$-harmonic
function $u$ admits an expansion in spherical harmonics

$$u(r\zeta)=\sum_{l\geq0}f_l(r^2)r^lu_l(\zeta)\eqno(3.6)$$
where $u_l$ is a spherical harmonic of degree $l$ and $f_l$ is the
hypergeometric function

$$f_l(x)=\frac{{}_2F_l(l,1-\frac{n}{2},l+\frac{n}{2},x)}{{}_2F_l
(l,1-\frac{n}{2},l+\frac{n}{2},1)}=\sum_{k=0}^\infty
\frac{\Gamma(l+k)\Gamma(1-\frac{n}{2}+k)
\Gamma(l+\frac{n}{2})\Gamma(1)}{\Gamma(l)\Gamma(1-\frac{n}{2})
\Gamma(l+\frac{n}{2}+k)\Gamma(1+k)}x^k.$$
Moreover the sum (3.6) converges uniformly on compact
subsets of $\B_n$, in particular

$$\norm{u_l}_{L^2(\S{n-1})}f_l(r^2)r^l=
\int_{\S^{n-1}}u(r\zeta)u_l(\zeta)d\sigma(\zeta).$$
On the other hand, if $l\not=0$ as $n$ is odd,

$$\frac{\Gamma(l+k)\Gamma(1-\frac{n}{2}+k)
\Gamma(l+\frac{n}{2})\Gamma(1)}{\Gamma(l)\Gamma(1-\frac{n}{2})
\Gamma(l+\frac{n}{2}+k)\Gamma(1+k)}=
\frac{\Gamma(l+\frac{n}{2})\Gamma(1)}{\Gamma(l)\Gamma(1-\frac{n}{2})}
\frac{1}{k^n}\ent{1+O\left(\frac{1}{k}\right)},$$
thus the $n-2$ first derivatives of $F_l$ have a limit when
$x\rightarrow1$,
whereas the $n-1$-st derivative grows like
$\log(1-x)$ when $x\rightarrow1$, thus $(3)$ implies that $u_l=0$ for
$l\not=0$, that is $u$ is constant.\hfill$\diamond$

$\diamond$ Assume now $n$ is even and write $n=2p$. Then if
$\ffi\in\cc^\infty(\S^{n-1})$, $\ffi$ admits a decomposition into spherical
harmonics
$\ffi=\sum_{l=0}^{+\infty}\ffi_l$ with $\norm{\ffi_l}_\infty=O(l^{-\alpha})$ for
every $\alpha>0$ (\cite{ST1} appendix C). But then

$$\P_h\ent{\ffi}(r\zeta)=\sum_{l=0}^{+\infty}f_l(r)r^l\ffi_l(\zeta)$$
with

$$f_l(r)r^l=\frac{{}_2F_1(l,1-p,l+p,r^2)}{{}_2F_1(l,1-p,l+p,1)}r^l=
\frac{\Gamma(l+2p-1)\Gamma(p)}{\Gamma(l+p)\Gamma(2p-1)}
\sum_{j=0}^p\frac{(l)_j(1-p)_j}{(l+p)_jj!}r^{2j+l}.$$
But, for every $k\geq 0$,

$$N^k\left(\sum_{j=0}^p\frac{(l)_j(1-p)_j}{(l+p)_jj!}r^{2j+l}\right)
=\sum_{j=0}^p\frac{(l)_j(1-p)_j}{(l+p)_jj!}(2j+l)^k2^kr^{2j+l}.$$
Therefore $N^k(f_lr^l)(1)=O(l^{k+p-1})$. But
$\norm{\ffi_l}_\infty=O(l^{-(k+p+1)})$ thus 
$\displaystyle\sum_{l=0}^{+\infty}N^kf_l(r)\ffi_l(\zeta)$
converges uniformly on $\overline{\B_n}$ and
$\P_h\ent{\ffi}\in\cc^\infty(\overline{\B_n})$.

The fact that for $u$ $\hh$-harmonic with a boundary distribution, $N^ku$ has
also a boundary distribution then results from the symmetry of the Poisson 
kernel : $\P_h(r\zeta,\xi)=\P_h(r\xi,\zeta)$.\hfill$\Box$
\end{proof}

\begin{remarquenum}[1] Normal derivatives of
$\hh$-harmonic functions have two opposite behaviors depending on the
dimension of $\B_n$. In odd dimension, the behavior is similar to the
complex case (see \cite{BBG}, in this case, the analog of constant
functions are pluriharmonic functions).

In opposite, in even dimension, the behavior is similar to that of
Euclidean harmonic functions.
\end{remarquenum}

\begin{remarquenum}[2] The similarity with the Euclidean case can be
seen in a different way. In \cite{Sam}, the following link between 
Euclidean harmonic functions and
$\hh$-harmonic functions has been proved :

\begin{lemmenum}[9] For every $\hh$-harmonic function $u$, there
exists
a unique Euclidean harmonic function $v$ such that $v(0)=0$ and :

$$u(r\zeta)=u(0)+\int_0^1v(rt\zeta)
\ent{(1-t)(1-tr^2)}^{\frac{n}{2}-1}\frac{dt}{t}$$
for every $0\leq r<1$ and every $\zeta\in\S^{n-1}$.
\end{lemmenum}

Moreover, let $f=\sum_lu_l$ is the spherical harmonics
expansion of $f\in L^2(\S^{n-1})$ and if
$g=\sum_l\frac{\Gamma(l+n-1)}{\Gamma(n-1)\Gamma(l)}u_l$, then lemma 9 
links  $u=\P_h\ent{f}$ to $v=\P_e\ent{g}$.

But, if $f=\sum_lu_l\in\cc^\infty(\S^{n-1})$ and 
$g=\sum_l\frac{\Gamma(l+n-1)}{\Gamma(n-1)\Gamma(l)}u_l$. Then, as
$\norm{u_l}_\infty=O(l^{-\alpha})$ for every $\alpha>0$,
$g\in\cc^\infty(\S^{n-1})$ thus
$v=\P_e\ent{g}\in\cc^\infty(\overline{\B_n})$.

Moreover,
if $n$ is even $(1-tr^2)^{\frac{n}{2}-1}$ is a polynomial and is
therefore
$\cc^\infty$, we then find again that $u\in\cc^\infty(\overline{\B_n})$.

In opposite, if $n$ is odd, we find again the $n-1$ obstacle since the
highest
order term of $(1-t)^{\frac{n}{2}-1}N^k(1-tr^2)^{\frac{n}{2}-1}$ is

$$(1-t)^{\frac{n}{2}-1}(1-tr^2)^{\frac{n}{2}-1-k}\simeq(1-t)^{n-2-k}$$
when $r\rightarrow1$, and since $(1-t)^{n-2-k}$ is not integrable for
$k\geq n-1$.
\end{remarquenum}

\section{Atomic decomposition of $\hh^p$ spaces}
\label{decompatom}

{\sl In this section we prove that $\hh^p$ spaces admit an atomic 
decomposition. In \ref{HpatHp} we define $\hh^p_{at}$ and show that this 
space is included in $\hh^p$. Conversely, we have seen in the previous 
chapter that $\hh$-harmonic functions in $\hh^p$ are obtained by $\hh$-Poisson integration
of distributions on $\S^{n-1}$, thus they are extensions of distributions
from $\S^{n-1}$ to $\B_n$. An other mean to extend a distribution on 
$\S^{n-1}$ to $\B_n$ is integration with respect to the Euclidean Poisson
kernel. In \ref{hareucharhyp} we study the links between this two extensions,
which allows us in \ref{hphat} to obtain the inclusion
$\hh^p\subset\hh^p_{at}$ from the atomic decomposition of $H^p$ spaces of 
Euclidean harmonic functions.}

\subsection{Links between Euclidean harmonic functions and
$\hh$-harmonic functions}
\label{hareucharhyp}

We will now prove a ``converse'' to lemma 9. 

\begin{lemmenum}[10] There exists a function
$\eta:\ent{0,1}\times\ent{0,1}\mapsto\R^+$ such that
\begin{description}
\item[i] $\P_e(r\zeta,\xi)=\int_0^1\eta(r,\rho)\P_h(\rho r\zeta,\xi)d\rho$,

\item[ii] there exists a constant $C$ such that for every
$r\in\ent{0,1}$, $\int_0^1\eta(r,\rho)d\rho\leq C$.
\end{description}
\end{lemmenum}

\begin{proof}[Proof] Note that $\frac{1}{(x+y)^{\frac{n}{2}}}=c_n
\int_0^\infty\frac{z^{\frac{n}{2}-2}}{(x+y+z)^{n-1}}dz$.
Writing $X=2(1-<\zeta,\xi>)$, with an obvious abuse of language, we then
get

\begin{align}
\P_e(r,X)=&\frac{1-r^2}{\bigl((1-r)^2+rX\bigr)^{\frac{n}{2}}}
=\frac{1-r^2}{r^{\frac{n}{2}}}
\frac{1}{\ent{\frac{(1-r)^2}{r}+X}^{\frac{n}{2}}}\notag\\
=&\frac{1-r^2}{r^{\frac{n}{2}}}c_n\int_0^\infty
\frac{z^{\frac{n}{2}-2}}{\ent{X+\frac{(1-r)^2}{r}+z}^{n-1}}dz\notag
\end{align}
The following change of variable
$z=\frac{(1-\rho)^2}{\rho}-\frac{(1-r)^2}{r}=\frac{(r-\rho)(1-\rho
r)}{\rho r}$, leads to

\begin{align}
\P_e(r,X)=&\frac{1-r^2}{r^{\frac{n}{2}}}c_n\int_0^r
\frac{\ent{(r-\rho)(1-\rho r)}^{\frac{n}{2}-2}}{\ent{X+
\frac{(1-\rho)^2}{\rho}}^{n-1}
(\rho r)^{\frac{n}{2}-2}}\frac{1-\rho^2}{\rho^2}d\rho\notag\\
=&\frac{1-r^2}{r^{n-2}}c_n\int_0^r
\frac{\ent{(r-\rho)(1-\rho r)}^{\frac{n}{2}-2}(1-\rho^2)}{\ent{\rho X+
(1-\rho^2)}^{n-1}\rho^{1-\frac{n}{2}}}d\rho\notag\\
=&\frac{1-r^2}{r^{n-2}}c_n\int_0^r \P_h(\rho,X)(1-\rho^2)^{2-n}
\ent{(r-\rho)(1-\rho
r)}^{\frac{n}{2}-2}\rho^{\frac{n}{2}-1}d\rho\notag\\
=&c_n(1-r^2)\int_0^1\P_h(rs,X)(1-r^2s^2)^{2-n}
\ent{(1-s)(1-sr^2)}^{\frac{n}{2}-2}s^{\frac{n}{2}-1}ds\notag
\end{align}
We thus obtain $i/$ with

$$\eta(r,s)=c_n(1-r^2)(1-r^2s^2)^{2-n}
\ent{(1-s)(1-sr^2)}^{\frac{n}{2}-2}s^{\frac{n}{2}-1}.$$
Of course $\eta\geq0$ and one easily checks that $\int_0^1\eta(r,s)ds\leq C$,
since $n\geq3$.\hfill$\Box$
\end{proof}

\begin{corollaire}[11] Let $\eta$ be the function defined by lemma 10.
Let $f$ be a distribution on $\S^{n-1}$ and let $u=\P_h\ent{f}$ and
$v=\P_e\ent{f}$. Then $u$ and $v$ are linked by

$$v(r\zeta)=\int_0^1\eta(r,s)u(rs\zeta)ds.$$
In particular, if $u\in\hh^p$, then $v\in H^p(\B_n)$ and
$\norm{v}_{H^p(\B_n)}\leq
C\norm{u}_{\hh^p}$.
\end{corollaire}

\subsection{The inclusion $\hh^p_{at}\subset \hh^p$}
\label{HpatHp}
%\index{patom@$p$-atome! sur $\S^{n-1}$|(}
\begin{definition} A function $a$ on $\S^{n-1}$ is called a {\rm
$p$-atom}
on $\S^{n-1}$ if either $a$ is a constant or $a$ is supported in a ball 
$\tilde B(\xi_0,r_0)$ and if
\begin{description}
\item[$1$] $\abs{a(\xi)}\leq\sigma\ent{\tilde
B(\xi_0,r_0)}^{-\frac{1}{p}}$,
for almost every $\xi\in\S^{n-1}$,

\item[$2$] for every function $\Phi\in\cc^{\infty}(\S^{n-1})$,

$$\abs{\int_{\S^{n-1}}a(\xi)\Phi(\xi)d\sigma(\xi)}\leq
\norm{\nabla^{k(p)}\Phi}_{L^\infty\bigl(\tilde B(\xi_0,r_0)\bigr)}
r_0^{k(p)}\sigma
\ent{\tilde B(\xi_0,r_0)}^{1-\frac{1}{p}}$$
with $k(p)$ an integer strictly bigger than
$(n-1)\left(\frac{1}{p}-1\right)$.
\end{description}
\end{definition}
%\index{patom@$p$-atome! sur $\S^{n-1}$|)}

\begin{proposition}[12] There exists a constant $C_p$ such that, for
every $p$-atom $a$ on $\S^{n-1}$, $A=\P_h\ent{a}$ satisfies

$$\norm{A}_{\hh^p(\B_n)}\leq C_p.$$
\end{proposition}

\begin{proof}[Proof] Let $a$ be a $p$-atom on $\S^{n-1}$, with support
in $\tilde B(\xi_0,r_0)$. We want to estimate

\begin{align}
\int_{\S^{n-1}}\sup_{t\in\ent{0,1}}&\abs{\int_{\tilde B(\xi_0,r_0)}
\P_h(t\zeta,\xi)a(\xi)d\sigma(\xi)}^pd\sigma(\zeta)\notag\\
=&\int_{\tilde B(\xi_0,cr_0)}\sup_{t\in\ent{0,1}}\abs{\int_{\tilde
B(\xi_0,r_0)}
\P_h(t\zeta,\xi)a(\xi)d\sigma(\xi)}^pd\sigma(\zeta)\notag\\
&+\int_{\S^{n-1}\setminus\tilde B(\xi_0,cr_0)}\sup_{t\in\ent{0,1}}
\abs{\int_{\tilde B(\xi_0,r_0)}
\P_h(t\zeta,\xi)a(\xi)d\sigma(\xi)}^pd\sigma(\zeta)\notag\\
=&I_1+I_2\notag
\end{align}
with $c>1$ a constant. But, by H\"older's inequality,

\begin{align}
I_1=&\int_{\tilde
B(\xi_0,cr_0)}\sup_{t\in\ent{0,1}}\abs{\P_h\ent{a}(t\zeta)}^p
d\sigma(\zeta)
\leq c\sigma\bigl(\tilde B(\xi_0,cr_0)\bigr)^{1-\frac{p}{2}}
\ent{\int_{\tilde
B(\xi_0,cr_0)}\sup_{t\in\ent{0,1}}\abs{\P_h\ent{a}(t\zeta)}^2
d\sigma(\zeta)}^{\frac{p}{2}}\notag\\
\leq&c\sigma\bigl(\tilde B(\xi_0,cr_0)\bigr)^{1-\frac{p}{2}}
\norm{\P_h\ent{a}}^p_{\hh^2(\B_n)}
\leq c\sigma\bigl(\tilde B(\xi_0,cr_0)\bigr)^{1-\frac{p}{2}}
\norm{a}^p_{L^2(\S^{n-1})}\notag
\end{align}
since $\P_h$ is bounded $L^2(\S^{n-1})\mapsto \hh^2(\B_n)$. Using
property
$(1)$ of atoms, we see that

$$I_1\leq C\left(\frac{\sigma\bigl(\tilde
B(\xi_0,cr_0)\bigr)}{\sigma\bigl(\tilde
B(\xi_0,r_0)\bigr)}\right)^{1-\frac{p}{2}}\leq C_p.$$

Let us now estimate $I_2$. Using property $(2)$ of atoms, we have, for
$\zeta\in\S^{n-1}\setminus\tilde B(\xi_0,cr_0)$

\begin{align}
\abs{\int_{\tilde
B(\xi_0,r_0)}\P_h(t\zeta,\xi)a(\xi)d\sigma(\xi)}^p\leq&
r_0^{pk(p)}\norm{\nabla_\xi^{k(p)}\P_h(t\zeta,\xi)}^p_{L^{\infty}}\sigma\bigl(\tilde
B(\xi_0,r_0)\bigr)^{p-1}\notag\\
\leq&C_pr_0^{pk(p)}(1-t^2)^{n-1}\times \sup_{\xi\in\tilde B(\xi_0,r_0)}
\frac{1}{d(\zeta,\xi)^{p(n+k(p)-1)}}\sigma\bigl(\tilde
B(\xi_0,r_0)\bigr)^{p-1}\notag
\end{align}
thus

\begin{align}
I_2\leq&C_pr_0^{pk(p)}\sigma\bigl(\tilde B(\xi_0,r_0)\bigr)^{p-1}
\int_{\S^{n-1}\setminus\tilde B(\xi_0,cr_0)}
\sup_{\xi\in\tilde B(\xi_0,r_0)}\frac{1}{d(\zeta,\xi)^{p(n+k(p)-1)}}
d\sigma(\zeta)\notag\\
\leq&C_p\frac{r_0^{pk(p)}r_0^{(n-1)(p-1)}}{r_0^{\ent{p\left(1+\frac{k(p)}{n-1}\right)-1}(n-1)}}\notag
\end{align}
since $p(n+k(p)-1)>n-1$ {\it i.e.}
$k(p)>(n-1)\left(\frac{1}{p}-1\right)$. Thus $I_2\leq C_p$.
\hfill$\Box$
\end{proof}

\begin{remarquenum}[1] Condition $(2)$ implies with $\Phi=1$ that

$$\int_{\S^{n-1}}a(\xi)d\sigma(\xi)=0.$$
\end{remarquenum}

\begin{remarquenum}[2] Condition $(2)$ is equivalent to the {\it a
priori}
weaker condition :

\begin{description}
\item[$2'$] For every spherical harmonic $P$ of degree $\leq k(p)$,
$$\abs{\int_{\S^{n-1}}a(\xi)P(\xi)d\sigma(\xi)}\leq
\norm{\nabla^{k(p)}P}_{L^\infty\bigl(\tilde B(\xi_0,r_0)\bigr)}
r_0^{k(p)}\sigma
\ent{\tilde B(\xi_0,r_0)}^{1-\frac{1}{p}}.$$
\end{description}

\begin{proof}[Proof] Assume this condition is fulfilled and let
$\Phi\in\cc^\infty(\S^{n-1})$. There exists $P$, a linear combination of
spherical harmonics of degree $\leq k(p)$ and $R\in\cc^\infty(\S^{n-1})$
such that

\begin{enumerate}
\item $\Phi=P+R$,

\item $\norm{R}_{\L^\infty\bigl(\tilde B(\xi_0,r_0)\bigr)}\leq C_p
r_0^{k(p)}
\norm{\nabla^{k(p)}\Phi}_{L^\infty\bigl(\tilde B(\xi_0,r_0)\bigr)}$.
\end{enumerate}
Then

\begin{align}
\abs{\int_{\S^{n-1}}a(\xi)\Phi(\xi)d\sigma(\xi)}\leq&
\abs{\int_{\S^{n-1}}a(\xi)P(\xi)d\sigma(\xi)}+
\abs{\int_{\S^{n-1}}a(\xi)R(\xi)d\sigma(\xi)}\cr
\leq&C_p\norm{\nabla^{k(p)}P}_{L^\infty\bigl(\tilde B(\xi_0,r_0)\bigr)}
r_0^{k(p)}\ent{\sigma\bigl(\tilde B(\xi_0,r_0)\bigr)}^{1-\frac{1}{p}}\cr
&+\norm{a}_{L^\infty\bigl(\tilde B(\xi_0,r_0)\bigr)}
\norm{R}_{L^\infty\bigl(\tilde B(\xi_0,r_0)\bigr)}
\sigma\bigl(\tilde B(\xi_0,r_0)\bigr)\cr
\leq&C\norm{\nabla^{k(p)}\Phi}_{L^\infty\bigl(\tilde B(\xi_0,r_0)\bigr)}
r_0^{k(p)}\ent{\sigma\bigl(\tilde B(\xi_0,r_0)\bigr)}^{1-\frac{1}{p}}.\notag\cr
\end{align}
\end{proof}
We could also impose the following weaker condition

\begin{description}
\item[$3$] For every spherical harmonic $P$ of degree $\leq
k(p)$,
$$\abs{\int_{\S^{n-1}}a(\xi)P(\xi)d\sigma(\xi)}=0.$$
\end{description}
We would than obtain a stronger atomic decomposition theorem. However
this
version is sufficient for our needs. It is also more intrinsic, the 
estimates we impose are directly those that are needed in the proof and
finally
it allows us to stay near to the proof in \cite{KL1}.
\end{remarquenum}

%\index{patom@$p$-atome! sur $\B_n$}
\begin{definition} A function $A$ on $\B_n$ is called an {\rm $\hh^p$-atom}
on $\B_n$ if there exists a $p$-atom $a$ on $\S^{n-1}$ such that
$A=\P_h\ent{a}$.

We define $\hh^p_{at}(\B_n)$ as the space of distributions $u$ on $\B_n$
such that there exists :

\begin{enumerate}
\item a sequence of $\hh_p$-atoms $(A_j)_{j=1}^{\infty}$ on $\B_n$,

\item a sequence $(\lambda_j)_{j=1}^{\infty}\in\ell^p$ such that

$$u=\sum_{j=1}^{\infty}\lambda_jA_j,\eqno(4.1)$$
with uniform convergence on compact subsets of $\B_n$.
\end{enumerate}

We write

$$\norm{u}_{\hh^p_{at}}=\inf\left\{\left(
\sum_{i=1}^\infty\abs{\lambda_j}^p\right)^{\frac{1}{p}}\right\}$$
where the infimum is taken over all decompositions of $u$ of the form
(4.1).
\end{definition}

\begin{proposition}[13] For $0<p\leq 1$,
$\hh^p_{at}(\B_n)\subset\hh^p(\B_n)$ there exists a constant $C_p$
such that for every $u\in\hh^p_{at}(\B_n)$,

$$\norm{u}_{\hh^p}\leq C_p\norm{u}_{\hh^p_{at}}.$$
\end{proposition}

\begin{proof}[Proof] It is {\it mutatis mutandis} the proof of theorem
2.2 in \cite{KL1}.

Let $\eps>0$ and let $u=\sum_{j=1}^{\infty}\lambda_jA_j$ be a function
in
$\hh^p_{at}$ and take an atomic decomposition such that
$\sum_{i=1}^\infty\abs{\lambda_j}^p\leq(1+\eps)\norm{u}_{\hh^p_{at}}^p$.

Property $2$ of atoms implies that

\begin{align}
\abs{\nabla^kA_j(x)}=&\abs{\nabla^k\P_h\ent{a_j}}\leq 
\norm{\nabla_\xi^{k(p)}\nabla_x^k\P_h(x,.)}_{L^\infty\bigl(B(\xi_0,r_0)\bigr)}
r_0^{k(p)}\sigma\bigl(B(\xi_0,r_0)\bigr)^{1-\frac{1}{p}}\cr
\leq&\frac{C_{p,k}}{(1-\abs{x})^{k_{p,l}}}\notag\cr
\end{align}
the series $\sum_{j=1}^\infty\lambda_j\nabla^kA_j(x)$ converge
uniformly on every compact subset of $\B_n$, thus
$\sum_{j=1}^\infty\lambda_jA_j(x)$ defines an $\hh$-harmonic function on
$\B_n$. 

Moreover

$$\abs{\sum_{j=1}^\infty\lambda_jA_j(x)}^p\leq 
\sum_{j=1}^\infty\abs{\lambda_j}^p\abs{A_j(x)}^p.$$
Therefore

\begin{align}
\int_{\S^{n-1}}\sup_{0<r<1}\abs{\sum_{j=1}^\infty\lambda_jA_j(r\zeta)}^p
d\sigma(\zeta)\leq&\int_{\S^{n-1}}\sup_{0<r<1}
\sum_{j=1}^\infty\abs{\lambda_j}^p\abs{A_j(r\zeta)}^pd\sigma(\zeta)\notag\\
\leq&C_p^p\sum_{j=1}^\infty\abs{\lambda_j}^p\notag\\
\leq&(1+\eps)^pC_p^p\norm{u}_{\hh^p_{at}}^p\notag
\end{align}
which means that $\norm{u}_{\hh^p}\leq C\norm{u}_{\hh^p_{at}}^p$.
\hfill$\Box$
\end{proof}

\subsection{The inclusion $\hh^p\subset\hh^p_{at}$}
%\markboth{Fonctions $\hh$-harmoniques}{L'inclusion $\hh^p\subset \hh^p_{at}$}
\label{hphat}

We will here use the fact that the space $H^p(\B_n)$ of Euclidean
harmonic
functions $v$ such that $\mm\ent{v}\in L^p(\S^{n-1})$ admits an atomic
decomposition {\it i.e.} that for every function $v\in H^p$, there
exists a sequence $(\lambda_k)_{k\in\N}\in\ell^p$ and a sequence
$(a_k)_{k\in\N}$
of $p$-atoms on $\S^{n-1}$ such that

$$v(r\zeta)=\sum_{k\in\N}\lambda_k\P_e\ent{a_k}(r\zeta)\eqno(4.2)$$
and moreover

$$\norm{v}_{H^p}\simeq
\left(\sum_{k\in\N}\abs{\lambda_k}^p\right)^{\frac{1}{p}}.$$
This result is well known, however it seems difficult to find an
adequate reference. One may for instance adapt the proof of Garnett 
and Latter \cite{GL1} as outlined in \cite{Col}.

Let $u\in\hh^p$, then $u$ admits a boundary distribution $f$ and
$u=\P_h\ent{f}$. Then let $v=\P_e\ent{f}$. By lemma 10,
$v\in H^p(\B_n)$ thus $v$ admits an atomic decomposition {\it i.e.}
there exists a sequence $(\lambda_k)_{k\in\N}\in\ell^p$ and a sequence
$(a_k)_{k\in\N}$ of $p$-atoms on $\S^{n-1}$ such that $v$ is given by (4.2), 
thus

$$f=\sum_{k=0}^\infty\lambda_ka_k$$
in the sense of distribution. Therefore
$u=\P_h\ent{\sum\lambda_ka_k}=\sum\lambda_k\P_h\ent{a_k}$, the series
being convergent in $\hh^p$ by proposition 13. We have thus 
proved the following theorem :

\begin{theorem}[14] For every $0<p\leq1$, $\hh^p=\hh^p_{at}$ and the
norms are equivalent.
\end{theorem}

\bibliographystyle{plain}
\bibliography{harm}

\begin{thebibliography}{10}

\bibitem{ABC}
\textsc{Ahern P., Bruna, J. and Cascante C.}
\newblock {$H^p$-theory for generalized $\mm$-harmonic functions in the unit
  ball}.
\newblock {\em Indiana Univ. Math. J.}, 45:103--145, 1996.

\bibitem{BBG}
\textsc{Bonami A., Bruna, J. and Grellier S.}
\newblock {On Hardy, $BMO$ and Lipschitz spaces of invariant harmonic functions
  in the unit ball}.
\newblock {\em Proc. London Math. Soc.}, 77:665--696, 1998.

\bibitem{Col}
\textsc{Colzani L.}
\newblock {Hardy spaces on unit spheres}.
\newblock {\em Boll. U.M.I. Analisi Funzionale e Applicazioni VI}, IV -
  C:219--244, 1985.

\bibitem{Erd}
\textsc{Erd\'ely and al}, editor.
\newblock {\em {Higher Transcendental Functions I}}.
\newblock Mac Graw Hill, 1953.

\bibitem{FS1}
\textsc{Fefferman C. and Stein E.M.}
\newblock {$H^p$ spaces of several variables}.
\newblock {\em Acta Math.}, 129:137--193, 1972.

\bibitem{GL1}
\textsc{Garnett J.B. and Latter R.H.}
\newblock {The atomic decomposition for Hardy spaces in several complex
  variables}.
\newblock {\em Duke J. Math.}, 45:815--845, 1978.

\bibitem{Jathese}
\textsc{Jaming Ph}.
\newblock {\em {Trois probl\`emes d'analyse harmonique}}.
\newblock PhD thesis, Universit\'e d'Orl\'eans, 1998.

\bibitem{KL1}
\textsc{Krantz S.G. and Li S.Y.}
\newblock {On decomposition theorems for Hardy spaces on domains in $\C^n$ and
  applications}.
\newblock {\em J. Fourier Anal. and Appl.}, 2:68--107, 1995.

\bibitem{Lew}
\textsc{Lewis J.B.}
\newblock {Eigenfunctions on symmetric spaces with distribution-valued boundary
  forms}.
\newblock {\em Jour. Func. Anal}, 29:287--307, 1978.

\bibitem{Mi1}
\textsc{Minemura K.}
\newblock {Harmonic functions on real hyperbolic spaces}.
\newblock {\em Hiroshima Math. J.}, 3:121--151, 1973.

\bibitem{Mi2}
\textsc{Minemura K.}
\newblock {Eigenfunctions of the Laplacian on a real hyperbolic spaces}.
\newblock {\em J. Math. Soc. Japan}, 27:82--105, 1975.

\bibitem{Sam}
\textsc{Samii H.}
\newblock {\em {Les Transformations de Poisson dans la Boule Hyperbolique}}.
\newblock PhD thesis, Universit\'e Nancy 1, 1982.

\bibitem{ST1}
\textsc{Stein E.M.}
\newblock {\em {Singular Integrals and Differentiability Properties of
  Functions}}.
\newblock Princeton University Press, 1970.

\bibitem{vdBS}
\textsc{van den Ban E.P. and Schlichtkrull H.}
\newblock {Assymptotic expansions and boundary values of eigenfunctions on
  Riemannian symmetric spaces}.
\newblock {\em J. Reine angew. Math.}, 380:108--165, 1987.

\end{thebibliography}
\end{document}